\documentclass{elsart}

\usepackage{amssymb}

\newtheorem{theorem}{Theorem}[section]
\newtheorem{lemma}[theorem]{Lemma}
\newtheorem{corollary}[theorem]{Corollary}
\newtheorem{proposition}[theorem]{Proposition}
\newtheorem{definition}[theorem]{Definition}

\newtheorem{remark}[theorem]{Remark}

\begin{document}

\begin{frontmatter}



\title{On pairwise sensitive homeomorphisms}


\author{C. A. Morales\thanksref{label2}}
\ead{morales@impa.br}

\thanks[label2]{Partially supported by CNPq, FAPERJ and PRONEX/Dynam. Sys. from Brazil.}

\address{Instituto de Matem\'atica,
Universidade Federal do Rio de Janeiro,
P. O. Box 68530,
21945-970 Rio de Janeiro, Brazil.}

\begin{abstract}
We obtain properties of the pairwise sensitive homeomorphisms defined in \cite{cj}.
For instance, we prove that their sets of points with converging semi-orbits have measure zero, that
such homeomorphisms do not exist in a compact interval and, in the circle,
they are the Denjoy ones. Applications including alternative proofs of well-known facts in expansive systems
are given.
\end{abstract}

\begin{keyword}
Pairwise Sensitive \sep Expansive \sep Metric Space
\MSC Primary  37B05 \sep Secondary 28C15.
\end{keyword}
\end{frontmatter}

\section{Introduction}
\label{intro}

\noindent
Let $(X,d)$ a compact metric space and $\mu$ be a Borel probability measure of $X$
with induced product measure $\mu^{\otimes 2}$.
A homeomorphism $f: X\to X$ is {\em pairwise sensitive}, with respect to $\mu$,
if there is positive constant $\delta$ such that
for $\mu^{\otimes 2}$-a.e. $(x,y)\in X\times X$
there is $n\in \mathbb{Z}$ such that $d(f^n(x), f^n(y))\geq \delta$.
This definition was introduced in \cite{cj}
although in the more general context of non-invertible maps.
Sufficient conditions for pairwise sensitivity were then obtained
as, for instance, to be weakly mixing or ergodic with positive entropy
on measurable partitions formed by continuity sets (c.f. \cite{cj}).

In this paper we obtain some properties of the pairwise sensitivity homeomorphisms on compact metric spaces.
For instance, we prove that the set of points with converging semi-orbits
under these homeomorphisms has measure zero. Additionally,
there are no such homeomorphisms of compact intervals
and, in the circle, they are precisely the Denjoy ones.
These results allow us to obtain probabilistic proofs of some well-known results about expansive systems.

\section{Statements and proofs}

\noindent
Given a homeomorphism $f: X\to X$, $x\in X$ and $\delta>0$ we define the dynamic ball
$$
\Gamma_\delta(x)=\{y\in X:d(f^n(x),f^n(y))\leq \delta, \forall n\in \mathbb{Z}\}.
$$
(The notation $\Gamma^f_\delta(x)$ to indicate the dependence on $f$.)

We start with the following characterization of the pairwise sensitive homeomorphisms.

\begin{lemma}
\label{th1-new}
The following properties are equivalent for every homeomorphism $f$
of a metric space $X$:
\begin{enumerate}
\item[(a)]
$f$ is pairwise sensitive.
\item[(b)]
There is $\delta>0$ such that
\begin{equation}
\label{paratodo}
\mu(\Gamma_\delta(x))=0,\quad\quad
\forall x\in X.
\end{equation}
\item[(c)]
There is $\delta>0$ such that
\begin{equation}
\label{quasitodo}
\mu(\Gamma_\delta(x))=0,
\quad\quad
\forall \mu\mbox{-a.e. }
x\in X.
\end{equation}
\end{enumerate}
\end{lemma}

\begin{pf}
Given $\delta>0$ define
$\mathcal{A}_\delta$ and $f\times f:X\times X\to X\times X$ by
$$
\mathcal{A}_\delta=\{(x,y)\in X\times X:d(x,y)\leq\delta\}
\quad
\mbox{ and }
\quad
(f\times f)(x,y)=(f(x),f(y)).
$$
As noticed in \cite{cj}, $f$ is
pairwise sensitive if and only if there is
$\delta>0$ satisfying
\begin{equation}
\label{eq2-new-new}
\mu^{\otimes 2}
\left(
\bigcap_{n\in \mathbb{N}}(f\times f)^{-n}(\mathcal{A}_{\delta})
\right)=0.
\end{equation}
On the other hand, the following inequalities hold
$$
\bigcap_{n\in \mathbb{N}}(f\times f)^{-n}(\mathcal{A}_\delta)\subset
\bigcup_{x\in X}(\{x\}\times \Gamma_\delta(x))
\subset
\bigcap_{n\in \mathbb{N}}(f\times f)^{-n}(\mathcal{A}_{2\delta})
$$
so
$$
F_\delta(x,y)\leq \chi_{\Gamma_\delta(x)}(y)\leq F_{2\delta}(x,y),
$$
where $F_\delta$ and $\chi_C$ denotes the characteristic functions of
$\bigcap_{n\in \mathbb{N}}(f\times f)^{-n}(\mathcal{A}_\delta)$ and $C\subset X$ respectively.
Integrating the last expression we obtain
$$
\mu^{\otimes 2}
\left(
\bigcap_{n\in\mathbb{N}}(f\times f)^{-n}(\mathcal{A}_\delta)
\right)
\leq
\int_X\int_X  \chi_{\Gamma_\delta(x)}(y)   d\mu(y)d\mu(x)
$$
\begin{equation}
\label{eqqq2-new-new}
\leq
\mu^{\otimes 2}
\left(
\bigcap_{n\in\mathbb{N}}(f\times f)^{-n}(\mathcal{A}_{2\delta})
\right)
\end{equation}

Now suppose that (a) holds, i.e., $f$ is pairwise sensitive.
So, there is $\delta>0$ satisfying
(\ref{eq2-new-new}).
Then, the second inequality in
(\ref{eqqq2-new-new}) implies $\mu(\Gamma_{\frac{\delta}{2}}(x))=0$ for $\mu$-a.e. $x\in X$
whence (c) holds.
On the other hand, if (b) holds, i.e., if
there is $\delta>0$ satisfying (\ref{paratodo}),
then the first inequality in (\ref{eqqq2-new-new}) implies (\ref{eq2-new-new}) so $f$ is pairwise sensitive
proving (a).

Next we prove that (b) and (c) are equivalent.
Indeed, we only have to prove that (c) implies (b).
Suppose by contradiction that (c) holds but not (b).
Since (c) holds there is $\delta>0$ satisfying (\ref{quasitodo}), and, since (b) fails,
there is $x_0\in X$ such that $\mu(\Gamma_{\delta/2}(x_0))>0$.
Denote $X_\delta=\{x\in X:\mu(\Gamma_\delta(x))=0\}$
so $\mu(X_\delta)=1$.
Since $\mu$ is a probability we obtain
$X_\delta\cap \Gamma_{\frac{\delta}{2}}(x_0)\neq\emptyset$ so
there is $y_0\in \Gamma_{\frac{\delta}{2}}(x_0)$ such that $\mu(\Gamma_\delta(y_0))=0$.
Now if $x\in \Gamma_{\frac{\delta}{2}}(x_0)$ we have $d(f^i(x),f^i(x_0))\leq \frac{\delta}{2}$
(for all $i\in \mathbb{N}$) and, since
$y_0\in \Gamma_{\frac{\delta}{2}}(x_0)$, we obtain
$d(f^i(y_0),f^i(x_0))\leq \frac{\delta}{2}$ (for all $i\in \mathbb{N}$)
so
$d(f^i(x),f^i(y_0))\leq d(f^i(x),f^i(x_0))+d(f^i(x_0),f^i(y_0))\leq \frac{\delta}{2}+\frac{\delta}{2}=\delta$
(for all $i\in \mathbb{N}$)
proving $x\in \Gamma_{\frac{\delta}{2}}(y_0)$.
Therefore $\Gamma_{\frac{\delta}{2}}(x_0)\subset \Gamma_\delta(y_0)$
so $\mu(\Gamma_{\frac{\delta}{2}}(x_0))\leq \mu(\Gamma_\delta(y_0))=0$ which is absurd.
This proves that (b) and (c) are equivalent.
This proves the result.
\qed
\end{pf}

For the next remark we recall the classical definition of expansive homeomorphism \cite{u}:

\begin{definition}
A homeomorphism $f$ is {\em expansive} if there is a positive $\delta$
(called {\em expansivity constant}) such that
$x=y$ whenever $d(f^n(x),f^n(y))\leq \delta$ for all $n\in\mathbb{Z}$.
\end{definition}

\begin{remark}
\label{the-remark}
It follows from the definition that a homeomorphism is expansive if and only if
there is $\delta>0$ such that $\Gamma_\delta(x)=\{x\}$ for all $x\in X$.
This
suggests to define {\em measure-expansive homeomorphism}
as a homeomorphism $f$ for which there is a positive $\delta$ such that
$\mu(\Gamma_\delta(x))=0$ for all $x\in X$ (see the unpublished paper \cite{m'}).
However, Lemma \ref{th1-new} implies that this definition is in fact equivalent to pairwise sensitivity.
Nevertheless, it is still pertinent to compare the pairwise sensitive and expansive systems.
\end{remark}

The following definition is motivated by the definition of expansivity constant above.

\begin{definition}
\label{defff}
A {\em $\mu$-expansivity constant} of $f$ is
a constant $\delta$ satisfying either (b) or (c) in Lemma \ref{th1-new}.
\end{definition}

It follows that $f$ is pairwise sensitive with respect to $\mu$ if
and only if it has a $\mu$-expansivity constant. From this
we obtain the following corollary which can be also derived directly from the definition.

\begin{corollary}
 \label{expansive-coro}
Every expansive homeomorphism is pairwise sensitive with respect to any
nonatomic Borel probability measure.
\end{corollary}

Notice that the converse of this corollary is false. More precisely,
there are homeomorphisms on certain compact metric spaces which are pairwise sensitive with respect to any
nonatomic Borel probability measure in $X_0$ but not expansive.
As a first example we can mention the identity $f: X_0\to X_0$ in $X_0=\{\frac{1}{n}:n\in \mathbb{N}^+\}\cup \{0\}$
with the metric induced by $\mathbb{R}$. Indeed since there is no
nonatomic Borel probability measures in $X_0$, $f$ is automatically pairwise sensitive
with respect to every nonatomic Borel probability measure but it is not expansive
(the author is indebted to a referee for suggesting this example).
A second example in which $X_0$ is a Cantor set is given in \cite{m}.

The following is nothing but a pairwise sensitive version of a well-known
characterization of the expansive homeomorphisms.

\begin{proposition}
\label{pp2}
If $n\in \mathbb{Z}\setminus\{0\}$, then a homeomorphism
$f$ is pairwise sensitive if and only if $f^n$ is.
\end{proposition}

\begin{pf}
It follows trivially from the definition of pairwise sensitivity that if $f^n$ is pairwise sensitive,
then so is $f$.
Conversely, suppose that $f$ is pairwise sensitive, and so,
it has a $\mu$-expansivity constant $\delta$. Since $X$ is compact and $n$ is fixed we can
choose $0<\epsilon<\delta$ such that if $d(x,y)\leq\epsilon$, then $d(f^i(x),f^i(y))<\delta$
for all $-n\leq i\leq n$.
With this property one has
$\Gamma_\epsilon^{f^n}(x)\subset\Gamma_\delta^f(x)$ for all $x\in X$
yielding $\mu(\Gamma_\epsilon^{f^n}(x))=0$ for $\mu$-a.e. $x\in X$. Thus $f^n$ is pairwise sensitive 
by Lemma \ref{th1-new}.
\qed
\end{pf}

If $f$ is a homeomorphism of $X$, and
$z\in X$, the {\em alpha-limit set} of $z$, denoted by $\alpha(z)$, consists of those points $y\in
X$ such that $y=\lim_{j\to\infty}f^{n_j}(z)$ for some sequence of integers
$n_j\to \infty$.
The {\em omega-limit set} of $z$, denoted by $\omega(z)$,
is similarly defined for
sequences $n_j\to-\infty$.
If, for some point $z$, $\alpha(z)$ and $\omega(z)$ each consist of a
single point, we say $z$ has {\em converging semi-orbits under $f$}.

As is well-known
if $f$ is an expansive homeomorphism of a compact metric space, then the set
of points having converging semi-orbits under $f$ is a countable set
(c.f. Theorem 1 in \cite{r} or Theorem 2.2.22 in \cite{ah}).
The following represents a pairwise sensitive version of this result.

\begin{theorem}
\label{reddy}
If $f$ is a pairwise sensitive homeomorphism with respect to $\mu$, then the set
of points with converging semi-orbits under $f$ has $\mu$-measure $0$.
\end{theorem}

A referee suggested
a proof of this theorem as in
Theorem 2.2.22 of the Aoki-Hiraide's book \cite{ah}.
Nevertheless, it seems that this does not apply in the present case
since the proof in \cite{ah} is based on the finiteness of the fixed point set.
Instead, we present the following alternative proof.

\begin{pf}
Denote by Fix$(f)$ and $A(f)$ the set of fixed points  and points with
converging semi-orbits under $f$ respectively.
Given $x,y\in X$, $n\in \mathbb{N}^+$ and $m\in \mathbb{N}$ we define
$A(x,y,n,m)$ as the set of points $z\in X$ satisfying
$$
\max\left\{d(f^i(z),x),d(f^j(z),y)\}\leq \frac{1}{n},\quad\forall i\leq-m\leq m\leq j\right\}.
$$
This set is clearly compact.

Let us show that there is a sequence of fixed points $x_k$ satisfying
\begin{equation}
\label{eq-reddy}
A(f)\subset\bigcap_{n\in \mathbb{N}^+}\bigcup_{k,k',m\in \mathbb{N}^+}A(x_k,x_{k'},n).
\end{equation}
Indeed, Fix$(f)$ is compact since $f$ is continuous so there is a sequence of fixed points $x_k$
which is dense in Fix$(f)$ (we shall prove that this sequence works).
Take $z\in A(f)$. Then, there are fixed points $x,y$ such that $\alpha(z)=x$ and $\omega(z)=y$.
Fix $n\in \mathbb{N}^+$. As $\alpha(z)=x$ and $\omega(z)=y$ there is
$m\in \mathbb{N}^+$ such that
$$
\max\{d(f^i(z),x),d(f^j(z),y)\}\leq \frac{1}{2n}
$$
whenever $i\leq -m\leq m\leq j$.
But $x,y\in$ Fix$(f)$ so there are $k,k'\in \mathbb{N}^+$ such that
$$
\max\{d(x,x_k),d(y,x_{k'})\}\leq\frac{1}{2n}.
$$
Combining these inequalities with the triangle inequality we obtain
$$
\max\{d(f^i(z),x_k),d(f^j(z),x_{k'})\}\leq \frac{1}{n},\quad\quad\forall i\leq -m\leq m\leq j
$$
proving
$z\in A(x_k,x_{k'},n,m)$ and so (\ref{eq-reddy}) holds.

Now, assume by contradiction that $\mu(A(f))>0$.
Then, (\ref{eq-reddy}) implies
\begin{equation}
\label{eq-reddy2}
\mu\left(\bigcup_{k,k',m\in \mathbb{N}^+}A(x_k,x_{k'},n,m)\right)>0,
\quad\quad\forall n\in \mathbb{N}^+.
\end{equation}
Since $f$ is pairwise sensitive we can fix a $\mu$-expansivity constant $e$ (c.f. Definition \ref{defff}).
Fix also $n\in \mathbb{N}^+$ such that $\frac{1}{n}\leq \frac{e}{2}$.
Applying (\ref{eq-reddy2}) to this $n$ we can arrange
$k,k',m\in \mathbb{N}^+$ such that
$$
\mu(A(x,y,n,m))>0,
$$
where $x=x_k$ and $y=x_{k'}$.
Then, since
$A(x,y,n,m)$ is compact,
there are $z\in A(x,y,n,m)$ and $\delta_0>0$ such that
$$
\mu(A(x,y,n,m)\cap B[z,\delta])>0,
\quad\quad\forall 0<\delta<\delta_0,
$$
where $B[\cdot,\delta]$ indicates the closed $\delta$-ball operation.
As $f$ is continuous and $m$ fixed
we can also fix $0<\delta<\delta_0$ such that
$$
d(f^i(z),f^i(w))\leq \frac{e}{2}
\mbox{ whenever } -m\leq i\leq m
\mbox{ and } d(z,w)<\delta.
$$
We claim that
$$
A(x,y,n,m)\cap  B[z,\delta]\subset \Gamma_e(z).
$$
Indeed, take $w\in A(x,y,n,m)\cap B[z,\delta]$.

Since $w\in B[z,\delta]$ one has $d(z,w)<\delta$ so
$$
d(f^i(w),f^i(z))\leq e,
\quad\forall -m\leq i\leq m.
$$

Since $z,w\in A(x,y,n,m)$
and $\frac{1}{n}\leq\frac{e}{2}$ one has
$$
d(f^i(w),f^i(z))\leq d(f^i(w),x)+d(f^i(z),x)\leq e
$$
and
$$
d(f^j(w),f^j(z))\leq d(f^j(w),y)+d(f^j(z),y)\leq e
$$
for all $i\leq -m\leq m\leq j$ proving $w\in \Gamma_e(z)$.
This proves the claim.

The claim implies
$$
0<\mu(A(x,y,n,m)\cap B[z,\delta])\leq \mu(\Gamma_e(z))
$$
which is absurd by Lemma \ref{th1-new} since $e$ is a $\mu$-expansivity constant.
This ends the proof.
\qed
\end{pf}

A direct corollary is the following pairwise sensitive version of Theorem 3.1 in \cite{u}.
Denote by Per$(f)$ the set of periodic points of $f$.

\begin{corollary}
\label{thA}
Every pairwise sensitive homeomorphism $f$ is {\em aperiodic}, i.e.,
$\mu($Per$(f))=0$.
\end{corollary}
 
\begin{pf}
Recalling Fix$(f)=\{x\in X:f(x)=x\}$ we have
Per$(f)=\cup_{n\in \mathbb{N}^+}$Fix$(f^n)$.
Now, $f^n$ is pairwise sensitive by Proposition \ref{pp2},
and every element of Fix$(f^n)$ is a point with converging semi-orbits of $f^n$,
thus $\mu($Fix$(f^n))=0$ for all $n$ by Theorem \ref{reddy}.
Therefore,
$\mu($Per$(f))\leq \sum_{n\in \mathbb{N}^+}\mu($Fix$(f^n))=0$.
\qed
\end{pf}

Now, we describe pairwise sensitive homeomorphisms in dimension one.
To start with we prove that there are no such homeomorphisms of compact intervals.

\begin{theorem}
\label{thD}
There are no pairwise sensitive homeomorphisms of a compact interval.
\end{theorem}

\begin{pf}
Suppose by contradiction that there is a
pairwise sensitive homeomorphisms of a compact interval $I$,
with respect to some Borel probability $\mu$ of $I$.
By Proposition \ref{pp2} we can assume that $f$ is orientation-preserving.
Since $f$ is also continuous we have that Fix$(f)\neq \emptyset$.
But such a set is also closed since $f$ is continuous so
its complement $I\setminus$Fix$(f)$ consists of countably many open intervals $J$.
Since $f$ is orientation-preserving it is clear that every point in $J$ is a point with converging
semi-orbits therefore $\mu(I\setminus$ Fix$(f))=0$ by Theorem \ref{reddy}.
But $\mu($Fix$(f))=0$ by Corollary \ref{thA} so
$\mu(I)=\mu($Fix$(f))+\mu(I\setminus$ Fix$(f))=0$ contradiction.
This ends the proof.
\qed
\end{pf}

\vspace{5pt}

Let $S^1$ denote the unit circle.
Recall that an orientation-preserving homeomorphism of the circle $S^1$ is {\em Denjoy}
if it is not topologically conjugated to a rotation \cite{hk}.

\begin{theorem}
\label{circle1}
A homeomorphism of $S^1$ is pairwise sensitive
(with respect to some Borel probability measure $\mu$) if and only if it is Denjoy.
\end{theorem}

\begin{pf}
Let $f$ be a Denjoy homeomorphism of $S^1$.
As is well known $f$ has no periodic points and exhibits a unique minimal set
$\Delta$ which is a Cantor set \cite{hk}.
In particular, $\Delta$ is compact without isolated points thus it exhibits nonatomic
Borel probability measures $\mu$ (c.f. Corollary 6.1 in \cite{prv}).
On the other hand, one sees that $f/\Delta$ is
an expansive homeomorphism, and so,
it is pairwise sensitive with respect to $\mu$ by Corollary \ref{expansive-coro}.

Conversely, suppose that $f$ is a pairwise sensitive homeomorphism of $S^1$
with respect to some Borel probability measure $\mu$.
Suppose by contradiction that it is not Denjoy.
Then, either $f$ has periodic points or is conjugated to a rotation (c.f. \cite{hk}).

In the first case we can assume by Proposition \ref{pp2} that
$f$ has a fixed point.
Then, cutting open $S^1$ along the fixed point we obtain a pairwise sensitive
homeomorphism of a compact interval
$I$ with respect to some Borel probability measure $\nu$.
But this contradicts Theorem \ref{thD}.

In the second case it is well-known that there are a homeomorphism $\phi$ and a
rotation $R$ of $S^1$ such that $\phi\circ f=R\circ \phi$ (c.f. \cite{hk}).
Let $\phi_*(\mu)$ be the pullback measure of $S^1$ defined by
$\phi_*(\mu)(A)=\mu(\phi^{-1}(A))$ for all borelian $A$.
Since $f$ is pairwise sensitive with respect to $\mu$ we have that $R$ is
pairwise sensitive with respect to $\phi_*(\mu)$.
However, isometries cannot be pairwise sensitive with respect to any
Borel probability measure (c.f. the remark after Theorem 2.2 in \cite{cj}).
This contradiction proves the result.
\qed
\end{pf}

\vspace{5pt}

We finish with some short applications of our results.
First we observe that combining Corollary 6.1 p. 210 in \cite{prv} and
Corollary \ref{expansive-coro} we obtain that
the set of periodic points of an expansive homeomorphism of a compact metric space is
countable (this was originally proved by Utz \cite{u}).
Second we obtain a probabilistic proof of
a result by Jacobsen and Utz (c.f. \cite{bry}, \cite{ju}).

\begin{corollary}
There are no expansive homeomorphisms of a compact interval.
\end{corollary}

\begin{pf}
Suppose by contradiction that there is an expansive homeomorphism of a compact interval $I$.
Since the Lebesgue measure $Leb$ is nonatomic we obtain that $f$ is
pairwise sensitive with respect to $Leb$ by Corollary \ref{expansive-coro}.
But such homeomorphisms do not exist by Theorem \ref{thD}.
\qed
\end{pf}

Recall that the {\em nonwandering set} of $f$ is the set $\Omega(f)$ consisting of those points
$x\in X$ for which $U\cap (\cup_{n\in \mathbb{N}^+}f^n(U))\neq\emptyset$ for all neighborhood $U$ of $x$.
As is well known $\Omega(f)$ contains the support supp$(\mu)$ of every invariant Borel probability measure
$\mu$ of $f$.
Now we prove the same but for pairwise sensitive homeomorphisms
in $S^1$.

\begin{lemma}
\label{supp}
If $f$ is a homeomorphism of $S^1$ which is pairwise sensitive with respect
to a Borel probability $\mu$ of $S^1$, then supp$(\mu)\subset\Omega(f)$.
\end{lemma}

\begin{pf}
Suppose by contradiction that
there is $x\in$ supp$(\mu)\setminus \Omega(f)$.
By Lemma \ref{th1-new} we have that $f$ has a $\mu$-expansivity constant $\delta$.
Since $x\notin\Omega(f)$ we can assume that
the collection of open intervals $\{f^n(B(x,\delta)):n\in \mathbb{Z}\}$ is disjoint.
Therefore, there is $N\in \mathbb{N}$ such that the length of $f^n(B(x,\delta))$
is less than $\delta$ for $|n|\geq N$.
From this and the continuity of $f$ we can arrange $\epsilon>0$
such that $B(x,\epsilon)\subset \Gamma_\delta(x)$
therefore $\mu(\Gamma_\delta(x))\geq \mu(B(x,\epsilon))>0$ as
$x\in$ supp$(\mu)$.
This contradicts Lemma \ref{th1-new} and the result follows.
\qed
\end{pf}

We use this lemma together with Theorem \ref{circle1} to obtain a probabilistic  proof
of the following result also by Jacobsen and Utz \cite{ju}.

\begin{corollary}
There are no expansive homeomorphisms of $S^1$.
\end{corollary}

\begin{pf}
Suppose by contradiction that there is an expansive homeomorphism of $S^1$.
Since the Lebesgue measure $Leb$ of $S^1$ is nonatomic we obtain from Corollary \ref{expansive-coro} that $f$ is
pairwise expansive with respect to $Leb$ and so supp$(Leb)\subset \Omega(f)$ by Lemma \ref{supp}.
But, $\Omega(f)$ is a Cantor set by Theorem \ref{circle1} and supp$(Leb)=S^1$ thus
we obtain a contradiction. This ends the proof.
\qed
\end{pf}


\bibliographystyle{amsplain}

\end{document}